\documentclass{amsart}\usepackage{amsmath}

\usepackage{amsfonts}

\newtheorem{theorem}{Theorem}
\theoremstyle{plain}

\newtheorem{remark}{Remark}

\numberwithin{equation}{section}

\input{tcilatex}

\begin{document}
\title[Universal enveloping locally C*-algebra for a locally JB-algebra]{On
existence and uniqueness of universal enveloping locally C*-algebra for a
locally JB-algebra}
\author{Alexander A. Katz}
\address{Alexander A. Katz, Dep. of Math \& CS, St. John's College of LAS,
St. John's University, 300 Howard Ave., DaSilva 314, Staten Island, NY
10301, USA}
\email{katza@stjohns.edu}
\author{Oleg Friedman}
\address{Oleg Friedman, Dep. of Math Sciences, University of South Africa,
P.O. Box 392, Pretoria 0003, South Africa}
\email{friedman001@yahoo.com}
\curraddr{Oleg Friedman, Dep. of Math \& CS, St. John's College of LAS, St.
John's University, 8000 Utopia Prkwy., Queens 11439, USA}
\email{fridmano@stjohns.edu}
\date{January 4, 2008}

\begin{abstract}
A theorem is presented on existence and uniqueness up to the topological
*-isomorphism of universal locally C*-algebra for an arbitrary locally
JB-algebra. 
\end{abstract}

\maketitle

The abstract Banach associative symmetrical *-algebras over $%
\mathbb{C}
,$ so called C*-algebras, were introduces first by Gelfand and Naimark in [$%
^{1}$]. In the present time the theory of C*-algebras become a vast portion
of Functional Analysis having connections and applications in almost all
branches of Modern Mathematics and Theoretical Physics (see for example [$%
^{2}$] for the basic theory of C*-algebras).

From the 1940's and the beginning of 1950's there were numerous attempts
made to extend the theory of C*-algebras to a category wider than Banach
algebras. For example, in 1952, while working on the theory of
locally-multiplicatively-convex algebras as projective limits of projective
families of Banach algebras, Arens in the paper [$^{3}$] and Michael in the
monograph [$^{4}$] independently for the first time studied projective
limits of projective families of functional algebras in the commutative case
and projective limits of projective families of operator algebras in the
non-commutative case. In 1971 Inoue in the paper [$^{5}$] explicitly studied
topological *-algebras which are topologically *-isomorphic to projective
limits of projective families of C*-algebras and obtained their basic
properties. He as well suggested a name of \textit{locally C*-algebras} for
that category. Below we will denote these algebras as \textit{LC*-algebras}.
For the present state of the theory of LC*-algebras see recently published
monograph of Fragoulopoulou [$^{6}$].

At the same time there were numerous attempts to extend the theory of
C*-algebras to non-associative algebras which are close to associative, in
particular to Jordan algebras. In fact, in 1978 Alfsen, Schultz and St\o %
rmer published their celebrated paper [$^{7}$], in which they introduced and
studied real Jordan Banach formally real algebras called \textit{JB-algebras}%
, which are real non-associative analogues of C*-algebras, and obtained for
this category analogues of the results from aforementioned paper [$^{1}$] by
Gelfand and Naimark. The exposition of elementary theory of JB-algebras can
be found in the monograph [$^{8}$] by Hanche-Olsen and St\o rmer, published
in 1984. In particular, in this monograph there is the following theorem
which was for the first time proved in 1980 by Alfsen, Hanche-Olsen and
Schultz in the paper [$^{9}$].

\begin{theorem}[{Alfsen, Hanche-Olsen, Schultz [$^{9}$]}]
For an arbitrary JB-algebra A there exists a unique up to an isometric
*-isomorphism a C*-algebra $C_{u}^{\ast }(A)$ (the \textit{universal
enveloping} C*-algebra for the JB-algebra A), and a Jordan homomorphism $%
\psi _{A}:A\rightarrow C_{u}^{\ast }(A)_{sa}$ from A to the self-adjoint
part of $C_{u}^{\ast }(A),$ such that: 

1). $\psi _{A}(A)$ generates $C_{u}^{\ast }(A)$ as a C*-algebra; 

2). for any pair composed of a C*-algebra $\mathfrak{A}$ and a Jordan
homomorphism $\varphi :A\rightarrow \mathfrak{A}_{sa}$ from A into the
self-adjoint part of $\mathfrak{A}$, there exists a *-homomorphism $\widehat{%
\varphi }:C_{u}^{\ast }(A)\rightarrow \mathfrak{A}$ from the C*-algebra $%
C_{u}^{\ast }(A)$ into C*-algebra $\mathfrak{A}$, such that $\varphi =%
\widehat{\varphi }\circ \psi _{A};$

3). there exists a *-antiautomorphism $\Phi $ of order 2 on the C*-algebra $%
C_{u}^{\ast }(A)$, such that $\Phi (\psi _{A}(a))=\psi _{A}(a),$ $\forall
a\in A.$
\end{theorem}

From the aforesaid one can see that it is natural and interesting to study
Jordan topological algebras which are projective limits of projective
families of JB-algebras. Those algebras under the name of \textit{locally
JB-algebras} were introduced and studied by Katz and Friedman in the paper [$%
^{10}$] published in 2006. In what follows we will call these algebras 
\textit{LJB-algebras}. 

An important question of the theory of LJB-algebras would have been an
analogue of the theorem 1 above. For the further exposition we need the
following technical theorem which is a corollary of Theorem 1 and general
properties of projective limits of projective families of Banach algebras.

\begin{theorem}
Let $\Lambda $ be a directed set of indices, and an arbitrary LJB-algebra A
be a projective limit $A=\underset{\longleftarrow }{\lim }A_{\alpha },$
where $\alpha \in \Lambda ,$ and $A_{\alpha }$ be a projective family of
JB-algebras. Then the family of C*-algebras $C_{u}^{\ast }(A)_{\alpha },$
where $\forall \alpha \in \Lambda ,$ $C_{u}^{\ast }(A)_{\alpha }$ be
universal enveloping C*-algebra for the corresponding JB-algebra $A_{\alpha
},$ is a projective family of C*-algebras.
\end{theorem}

Using the theorems 1 and 2 we are able to obtain the following main theorem
about existence and uniqueness of universal enveloping LC*-algebra for an
arbitrary LJB-algebra.

\begin{theorem}
For an arbitrary LJB-algebra A there exists a unique up to a topological
*-isomorphism a LC*-algebra $LC_{u}^{\ast }(A)$ (the \textit{universal
enveloping} locally C*-algebra for the LJB-algebra A), and a Jordan
homomorphism $\psi _{A}:A\rightarrow LC_{u}^{\ast }(A)_{sa}$ from A to the
self-adjoint part of $LC_{u}^{\ast }(A),$ such that: 

1). $\psi _{A}(A)$ generates $LC_{u}^{\ast }(A)$ as a LC*-algebra; 

2). for any pair composed of a LC*-algebra $\mathfrak{A}$ and a Jordan
homomorphism $\varphi :A\rightarrow \mathfrak{A}_{sa}$ from A into the
self-adjoint part of $\mathfrak{A}$, there exists a *-homomorphism $\widehat{%
\varphi }:LC_{u}^{\ast }(A)\rightarrow \mathfrak{A}$ from the LC*-algebra $%
LC_{u}^{\ast }(A)$ into LC*-algebra $\mathfrak{A}$, such that $\varphi =%
\widehat{\varphi }\circ \psi _{A};$

3). there exists a *-antiautomorphism $\Phi $ of order 2 on the LC*-algebra $%
LC_{u}^{\ast }(A)$, such that $\Phi (\psi _{A}(a))=\psi _{A}(a),$ $\forall
a\in A.$
\end{theorem}

\begin{remark}
The results of this short notes were presented by authors at the Fall
Meeting of the American Mathematical Society No. 1031 which took place on
October 6-7, 2007 in New Brunswick, NJ, USA.
\end{remark}

\begin{center}
REFERENCES
\end{center}

[$^{1}$] Gelfand, I.M.; Naimark, M.A., On the embedding of normed rings into
the ring of operators in Hilbert space. (English. Russian summary) Rec.
Math. [Mat. Sbornik] N.S. Vol. 12(54) (1943), pp. 197-213.

[$^{2}$] Pedersen, G.K., C*-algebras and their automorphism groups.
(English), London Mathematical Society Monographs. Vol. 14. London - New
York -San Francisco: Academic Press., 416 pp., (1979).

[$^{3}$] Arens, R., A generalization of normed rings. (English), Pac. J.
Math., Vol. 2 (1952), pp. 455-471.

[$^{4}$] Michael, E.A., Locally multiplicatively-convex topological
algebras. (English), Mem. Am. Math. Soc., Vol. 11 (1952), 79 pp. 

[$^{5}$].Inoue, A., Locally C*-algebras. (English), Mem. Fac. Sci. Kyushu
Univ. (Ser. A), No. 25 (1971), pp. 197-235.

[$^{6}$] Fragoulopoulou, M., Topological algebras with involution. (English)
North-Holland Mathematics Studies, Vol. 200. Elsevier Science B.V.,
Amsterdam, 495 pp., (2005).

[$^{7}$]. Alfsen, E.M.; Schultz, F.W.; St\o rmer, E., A Gelfand-Naimark
theorem for Jordan algebras. (English) Advances in Math. Vol. 28 (1978), No.
1, pp. 11-56.

[$^{8}$]. Hanche-Olsen, H.; St\o rmer, E., Jordan operator algebras.
(English), Monographs and Studies in Mathematics, Vol. 21. Boston - London -
Melbourne: Pitman Advanced Publishing Program. VIII, 183 pp., (1984).

[$^{9}$] Alfsen, E.M.; Hanche-Olsen, H.; Schultz, F.W., State spaces of
C*-algebras. (English) Acta Math. Vol. 144 (1980), No. 3-4, pp. 267-305. 

[$^{10}$] Katz, A.A.; Friedman, O., On projective limits of real C*- and
Jordan operator algebras. (English), Vladikavkaz Mathematical Journal, Vol.
8 (2006), No. 2, pp. 33-38.

\end{document}